\documentclass[a4paper,10pt]{article}
 \usepackage{hyperref}
\usepackage{paper-en}

\def\thetitle{Metric foundations of geometry}
\def\theauthors{Nina Lebedeva and Anton Petrunin}

\hypersetup{colorlinks=true,
citecolor=black,
linkcolor=black,
anchorcolor=black,
filecolor=black,
menucolor=black,
urlcolor=black,
pdftitle={\thetitle},
pdfauthor={\theauthors}
}

\begin{document}
\title{\thetitle}
\author{\theauthors}
\date{}
\maketitle

\begin{abstract}
A metric space is called all-set-homogeneous if every isometry between two of its subsets
extends to an isometry of the whole space.
We classify all-set-homogeneous length spaces: besides the classical examples, they include the universal metric trees of finite valence.
\end{abstract}

\section{Introduction}

A metric space is said to be \emph{all-set-homogeneous}
if any isometry between its subsets can be extended to an isometry of the whole space; here \emph{isometry} is a bijective distance-preserving map.

Examples of all-set-homogeneous spaces include complete simply connected Riemannian manifolds with constant curvature and circles with length metrics.
These will be referred to as \emph{classical spaces}.

The so-called \emph{universal metric trees} $\TT_n$ of finite valence $n\geqslant 2$ give another source of examples.
These are complete metric trees with valence $n$ at any point \cite{dyubina-polterovich}.
The space $\TT_n$ is uniquely defined up to isometry; it is homogeneous, and if $n$ is finite, then it is all-set-homogeneous; the last statement was proved in \cite{lebedeva-petrunin2024}.

All the \emph{terminology} used here, including metric trees, is explained in the next section.

\begin{thm}{Theorem}\label{thm:main}
Any all-set-homogeneous length space is a universal metric tree of finite valence or a classical space.
\end{thm}

Our theorem provides a very economical metric foundation for geometry.
A special case of our theorem was proved by Garrett Birkhoff \cite{birkhoff};
namely, he proved that \textit{any all-set-homogeneous geodesic space with locally unique nonbifurcating geodesics is classical}.

For locally compact spaces, our result follows from the structure theory of locally compact groups; moreover, \textit{any locally compact 3-point-homogeneous length space is classical};
a metric space is called \emph{$n$-point-homogeneous} if any isometry between its subsets with at most $n$ points in each can be extended to an isometry of the whole space.
This result follows from a theorem of Jacques Tits \cite{tits} about 2-point-homogeneous spaces;
under additional geometric assumptions, the 3-point-homogeneity statement was proved earlier by Herbert Busemann \cite[V §9]{busemann}.

Our only tools are symmetries, almost shortest paths, and repeated applications of the triangle inequality.
As in many classification problems, we must first prove a few curious facts about the hypothetical extra examples, only to discover that they form the empty set.

The theorem is reduced to the following statements:
\begin{enumerate}[(a)]
\item All-set-homogeneous length spaces are complete; see~\ref{prop:complete}.
\item Complete all-set-homogeneous length spaces are geodesic; see~\ref{prop:ash>geodesic}.
\item All-set-homogeneous geodesic spaces have locally unique geodesics; see~\ref{prop:unique-geod}.
\item Any all-set-homogeneous geodesic space with locally unique geodesics is a universal tree of finite valence or a classical space; see \ref{thm:tree}.
\end{enumerate}
These steps can be read independently; if the reader is interested only in uniquely geodesic spaces, then it is sufficient to read only the last pages of Section~\ref{Locally unique geodesics}.

The final section discusses an alternative approach to the proof and extensions beyond length spaces.
It can be read independently of the rest.

\section{Preliminaries}

This section should be used as a reference.

\subsection*{Terminology}

A metric space is called a \emph{length space} if the distance between any two points can be approximated arbitrarily well by lengths of curves joining them.

A \emph{geodesic} in a metric space $M$ is a distance-preserving map from a real interval to $M$.
The image of a geodesic with endpoints $x$ and $y$ can be denoted by $[xy]$ and will also be called a \emph{geodesic}.
If for any $x,y\in M$ there is a geodesic $[xy]$, then $M$ is called a \emph{geodesic space}.
If any point $p\in M$ admits a neighborhood $U$ such that the geodesic $[xy]$ is uniquely defined for any $x,y\in U$,
then $M$ has \emph{locally unique geodesics}.
If $[px]\subset [py]$ or $[px]\supset [py]$ whenever $[px]\cap [py]$ contains a nontrivial geodesic, then we say that geodesics in $M$ are \emph{nonbifurcating}.

Interior and closure in a metric space $M$ will be denoted by $\Int_M$ and $\Cl_M$, respectively.
A \emph{sphere}, an \emph{open ball}, and a \emph{closed ball} with center $p$ and radius $r>0$ in a metric space $M$ are defined as
\begin{align*}
\Sph(p,r)=\Sph(p,r)_M&=\set{x\in M}{|p-x|=r},
\\
\Ball(p,r)=\Ball(p,r)_M&=\set{x\in M}{|p-x|<r}
\\
\bar\Ball(p,r)=\bar\Ball(p,r)_M&=\set{x\in M}{|p-x|\leqslant  r};
\end{align*}
here $|\ -\ |$ denotes the distance between points in $M$.
We say that $x$ is a midpoint of $p$ and $q$ if $|p-x|=|q-x|=\tfrac12\cdot|p-q|$.

\subsection*{Observations}

The next three statements are essentially due to Piotr Niemiec \cite[3.9, 3.28, 3.31]{niemiec2023}.

\begin{thm}{Observation}\label{obs:spheres}
The intersection of an arbitrary family of spheres in an all-set-homogeneous metric space $M$ is all-set-homogeneous.
In particular, the set of midpoints of any pair of points in $M$ is all-set-homogeneous.
\end{thm}

\parit{Proof.}
Choose a family of spheres;
let $\Sigma$ be its intersection, and let $P$ be the set of their centers.

An isometry $X \to X'$ between subsets $X,X'\subset \Sigma$ can be extended by the identity on $P$.
It defines an isometry $P\cup X\to P\cup X'$,
which in turn extends to an isometry, say $\iota$, of $M$.
Since $\iota$ fixes the centers of the spheres, it defines an isometry of $\Sigma$.

The last statement follows since $\Sph(p,\tfrac12\cdot|p-q|)\cap \Sph(q,\tfrac12\cdot|p-q|)$ is the set of all midpoints of a pair $p,q\in M$.
\qeds

\begin{thm}{Observation}\label{obs:FIP}
Let $\Sigma_1,\Sigma_2,\dots$ be a sequence of spheres in an all-set-homogeneous metric space $M$.
If $\Sigma_1\cap\ldots\cap\Sigma_n\ne\emptyset$ for any $n$, then
\[
\Sigma_1\cap\Sigma_2\cap\ldots\ne\emptyset.
\]

\end{thm}

\parit{Proof.}
Suppose $\Sigma_i=\Sph(p_i,r_i)$; we can assume that $p_i\ne p_j$ if $i\ne j$.
Consider a countable set $X=\{x,y_1,y_2,\dots\}$ with metric defined by $|x-y_i|=r_i$ and $|y_i-y_j|=|p_i-p_j|_M$ for all $i$ and $j$.

Note that every finite subset $\{x,y_1,\ldots,y_n\}$ admits a distance-preserving map $f_n$ into $M$.
Since $M$ is all-set-homogeneous, we can assume $f_{n+1}$ is an extension of $f_n$ for every $n$.
Define $f\colon X\to M$ as the extension of all $f_n$.
Applying all-set-homogeneity, we can assume that $f(y_n)=p_n$ for every $n$.
Then $f(x)\in \Sigma_n$ for every $n$.
\qeds

\begin{thm}{Observation}\label{lem:length-sphere}
Assume an all-set-homogeneous metric space $X$ contains a nontrivial rectifiable curve. Then the following holds.
\begin{enumerate}[(i)]
\item There is $\eps>0$ such that any two connected components of $X$ lie $\eps$-apart from each other.
\item Let $Y$ be a connected component of $X$; equip it with the length metric $|\ \z-\ |_Y$ induced by the metric $|\ -\ |_X$ on $X$.
Then $Y$ is an all-set-homogeneous length space; moreover, there is a strictly increasing continuous function $\psi_X$ such that $\psi_X(0)=0$ and
\[|x-y|_Y=\psi_X(|x-y|_X)\]
for any $x,y\in Y$.
\end{enumerate}

\end{thm}

\parit{Proof.}
We can assume that the curve is compact; let $\eps>0$ be its diameter.

\parit{(i)} If $x$ and $y$ lie in different connected components of $X$, then $|x-y|_X>\eps$.

\parit{(ii)}  If $x,y\in Y$,
then there is a chain of points $x=x_0,x_1,\ldots,x_n=y$ such that $|x_i-x_{i+1}|\leqslant\eps$ for every $i$.
In this case, $x_i$ can be connected to $x_{i+1}$ by a rectifiable curve;
thus the same can be done for $x$ and $y$.
In particular, $|x-y|_Y<\infty$.

It remains to apply all-set-homogeneity of $X$.
\qeds

A subset $E$ of a metric space is called equilateral if the distance $|e-e'|$ is the same for any pair of distinct points $e,e'\in E$.

\begin{thm}{Observation}\label{obs:equilateral}
Every equilateral subset of an all-set-homogeneous metric space is finite.
\end{thm}

The proof repeats the argument of Garrett Birkhoff \cite[Theorem 18]{birkhoff}.
The same argument also implies that an all-set-homogeneous space cannot be isometric to a proper subset of itself; for example, an infinite-dimensional Hilbert space is not all-set-homogeneous.

\parit{Proof.}
Assume that an all-set-homogeneous metric space $M$ contains an infinite equilateral set $E$.
Extend $E$ to a maximal equilateral set $E_{\max}$ with respect to inclusion.
Choose a bijection $f$ from $E_{\max}$ onto a proper subset of it.
Note that $f$ is distance-preserving, and it does not extend to an isometry $F$ of $M$;
otherwise $F^{-1}(e)$ for $e\in E_{\max}\setminus f(E_{\max})$ can be added to $E_{\max}$, which is already maximal.
Therefore, $M$ is not all-set-homogeneous --- a contradiction.
\qeds

Let us define a \emph{local geodesic} to be a locally distance-preserving map from a real interval to a metric space.
A geodesic space $M$ will be called \emph{geodesically complete}
if any nontrivial geodesic in $M$ can be extended to a local geodesic defined on the whole real line.

\begin{thm}{Observation}\label{obs:geodesically-complete}
Any all-set-homogeneous geodesic space $M$ is geodesically complete.
\end{thm}

There are 2-point-homogeneous metric trees of any valence $n\geqslant 3$ that are not complete and not geodesically complete in the above sense, and hence 2-point-homogeneity is not sufficient.
These examples were constructed by Pénélope Azuelos~\cite{azuelos2025}, and they have the following weaker property: any geodesic defined on a closed interval can be extended to a geodesic defined on the whole real line.
The main theorem (once it is proved) will imply a stronger statement: \textit{every $(\diam M)$-long local geodesic in an all-set-homogeneous geodesic space $M$ is a geodesic}.

\parit{Proof.}
Choose a nontrivial geodesic $\gamma$ in $M$;
suppose it is defined on an interval~$\II$.

Choose a subinterval $\JJ\subset \II$ and its shift $\JJ'\subset \II$ by $a>0$ such that $\II=\JJ\cup\JJ'$.
The partial map $\gamma(t)\mapsto \gamma(t+a)$ extends to an isometry $F$ of $M$,
and we can extend $\gamma$ to a local geodesic defined by $\hat\gamma(t+n\cdot a)\df F^n\circ \gamma(t)$
for $t\in \JJ$ and every integer~$n$.
\qeds

\subsection*{Trees and cut points}

A uniquely geodesic space $T$ is a \emph{metric tree} if every geodesic triangle in $T$ is a \emph{tripod};
that is, the union of any two sides of any geodesic triangle in $T$ contains the remaining side.

\begin{thm}{Claim}\label{clm:complete}
Any geodesically complete metric tree $T$ is complete.
\end{thm}

\parit{Proof.}
Assume there is a nonconvergent Cauchy sequence $x_0,x_1,\ldots$ in~$T$.
Let $\ell=\lim |x_0-x_n|$;
note that $\ell>0$.

Denote by $\gamma_n$ the geodesic from $x_0$ to $x_n$ that is parametrized from $x_0$.
Since $T$ is a tree,
\[2\cdot t\leqslant|x_0-x_i|+|x_0-x_j|-|x_i-x_j|
\quad\Rightarrow\quad
\gamma_i(t)=\gamma_j(t).\]
Therefore the sequence $\gamma_n(t)$ stabilizes for any $0\leqslant t<\ell$.
Defining $\gamma_\infty(t)$ as $\gamma_n(t)$ for sufficiently large $n$ produces a geodesic $\gamma_\infty \colon[0,\ell)\to T$.
Since $T$ is geodesically complete, $\gamma_\infty$ extends to $[0,\ell]$, and then $x_0,x_1,\ldots$ converge to $\gamma_\infty(\ell)$ --- a contradiction.
\qeds

The \emph{valence} of a point $x$ in a metric tree $T$ is defined as the cardinality of the set of connected components in $T\setminus \{x\}$.
Recall that $\TT_n$ denotes the universal metric tree of valence $n$,
which is the (unique up to isometry) complete metric tree with valence $n$ at any point \cite{dyubina-polterovich}.

\begin{thm}{Proposition}\label{prop:global-tree}
Let $M$ be an all-set-homogeneous length space.
Suppose some (and therefore every) point of $M$ has a neighborhood isometric to a metric tree of valence at least $3$.
Then $M$ is a universal metric tree of finite valence.
\end{thm}

\parit{Proof.}
It is sufficient to show that $M$ does not contain an embedded circle; see for example \cite[2.2.3]{chiswell}.

Assume the contrary; let $2\cdot R$ be the infimum of lengths of embedded circles in $M$.
Since $M$ is locally a tree, homogeneity implies that $R>0$.
Absence of embedded circles of length less than $2\cdot R$ implies the following:
\begin{itemize}
\item an $\tfrac R2$-neighborhood of any point in $M$ is a tree with valence at least $3$;
\item if $|x-y|\z<R$, then $[xy]$ is uniquely defined.
\end{itemize}

It follows that any equilateral triangle in $M$ with side $s<R$ is a tripod.
Indeed, choose a point $p$ and three points $x$, $y$, and $z$ in different components of its punctured $\tfrac R2$-neighborhood at distance $\tfrac s2$ from $p$.
Note that $[xyz]$ is a tripod with side~$s$;
by $3$-point-homogeneity the same holds for any triangle with side $s$.

Now choose an embedded circle $\Sigma$ of length slightly above $2\cdot R$.
Choose points $x,y,z$ that divide $\Sigma$ into three equal arcs, say of length $s$;
we can assume  $s<R$.
Note that $|x-y|=s$; otherwise $\Sigma$ together with $[xy]$ would contain an embedded circle of length less than $2\cdot R$,
which is impossible.
It follows that $[xyz]$ is an equilateral triangle of side $s$,
which is not a tripod --- a contradiction.

Thus $M$ is an all-set-homogeneous metric tree;
in particular, it has the same valence at every point.
By \ref{obs:geodesically-complete} and \ref{clm:complete}, $M$ is complete.
Thus $M$ is a universal metric tree;
it remains to show that it has finite valence.

Assume the valence is infinite.
Then there is an infinite equilateral set $E\z\subset M$;
indeed, one can take points at a fixed distance in distinct branches of the tree.
Now we arrive at a contradiction to \ref{obs:equilateral}.
\qeds

Note that any open ball in a length space $M$ is path-connected.
A point $p\in M$ is called a \emph{local cut point} if some punctured ball $\Ball(p,R)\setminus\{p\}$ is not connected.
For us the empty set is not connected; so an isolated point of $M$ is a local cut point.

\begin{thm}{Proposition}\label{prop:cut-point}
If an all-set-homogeneous length space $M$ has a local cut point,
then it is isometric to a universal tree of finite valence,
a one-point space,
or a circle with a length metric.
In particular, $M$ is complete.
\end{thm}

\parit{Proof.}
Choose $p\in M$ and $r>0$ such that $\Ball(p,r)\setminus\{p\}$ is not connected.
We can assume that $\Ball(p,r)\setminus\{p\}\ne\emptyset$; otherwise $M$ is a one-point space.

Let $V$ and $W$ be two connected components of $\Ball(p,r)\setminus\{p\}$;
choose two points $v\in V$ and $w\in W$ and $0<\eps<\tfrac r2$ such that $|p-v|>\eps$ and $|p-w|>\eps$.
Since $M$ is a length space, we have
\[
0<s<\eps
\quad\Rightarrow\quad
\Sph(p,s)\cap V\ne\emptyset
\quad\text{and}\quad
\Sph(p,s)\cap W\ne\emptyset.
\]

Choose two points $x,y\in M$ such that $|x-y|<\eps$.
By all-set-homogeneity, given $0<s<|x-y|$, there is a unique point $z_s$ such that
\[|x-z_s|+|z_s-y|=|x-y|
\quad\text{and}\quad
|x-z_s|=s.\]
Moreover, $|z_{s'}-z_s|=|s'-s|$, hence there is a unique geodesic $[xy]$.

Now assume that a geodesic triangle $[xyz]$ with sides smaller than $\eps$ is not a tripod.
By all-set-homogeneity, we can assume that $x\in V$, $y\in W$, and
\[2\cdot |p-x|+|y-z|=|x-y|+|x-z|.\]
In this case $p\in [xy]$.
By the triangle inequality, we also have $p\in [yz]$ or $p\in [xz]$.
Hence $p\in [yz]$ and $p\in [xz]$ --- a contradiction.
It follows that an $\tfrac\eps2$-neighborhood of any point in $M$ is a metric tree.
Since $M$ is all-set-homogeneous, this tree has the same valence $n\geqslant 2$ at every point.

If $n=2$, then $M$ is locally isometric to the real line.
In this case $M$ must be isometric to the real line $\RR=\TT_2$ or a circle with a length metric.

If $n\geqslant3$, then by \ref{prop:global-tree}  $M$ is a universal metric tree of finite valence.
\qeds

\subsection*{Reflections and mirrors}

Let $M$ be an all-set-homogeneous space.
Define the \emph{mirror} of two distinct points $u$ and $v$ in $M$ as the set
\[\Bis(u,v)=\set{x\in M}{|u-x|=|v-x|}.\]

\begin{thm}{Observation}\label{lem:bisector}
Let $u$ and $v$ be distinct points in an all-set-homogeneous length space $M$.
Then $M\setminus \Bis(u,v)$ has two path-connected components that contain $u$ and $v$, respectively.
Furthermore, if the two values $|u-x|-|v-x|$ and $|u-y|-|v-y|$ have opposite signs or one of them vanishes, then
\[|x-y|=\inf\set{|x-o|+|o-y|}{o\in \Bis(u,v)}.\]
\end{thm}

The following observation is a version of Birkhoff's reflection construction
\cite[Theorem 20]{birkhoff}.

\begin{thm}{Observation}\label{lem:reflection}
Let $M$ be an all-set-homogeneous length space, and let $u,v\z\in M$ be distinct.
Then there is an involutive isometry $\rho$ that fixes $\Bis(u,v)$ pointwise and swaps $u$ with $v$.
Moreover,
\[\tfrac12\cdot |x-\rho(x)|=\inf\set{|x-o|}{o\in \Bis(u,v)}\]
for any $x\in M$.
Furthermore,
\[x\ne \rho(x)
\quad\Rightarrow\quad
\Bis(x,\rho(x))=\Bis(u,v).\]
\end{thm}

\parit{Proof.}
The partial isometry fixing $\Bis(u,v)$ and swapping $u$ with $v$ extends to an isometry, say $\iota$, of $M$.
It swaps the two sides of $\Bis(u,v)$ defined by the sign of $|u-x|-|v-x|$.
Define $\rho(x)=\iota(x)$ if $|u-x|\leqslant |v-x|$ and $\rho(x)=\iota^{-1}(x)$ otherwise.
By construction, $\rho$ is an involution and its fixed set is $\Bis(u,v)$.
Note that $\rho$ preserves distances to every point $o\in \Bis(u,v)$ and apply~\ref{lem:bisector}.

If $x\ne \rho(x)$, then repeating the construction above for $x$ and $\rho(x)$, we obtain an involutive isometry $\rho'$.
Note that $\Bis(u,v)\subset\Bis(x,\rho(x))$; therefore, $\rho'$ fixes $\Bis(u,v)$ pointwise.
The two strict sides of $\Bis(u,v)$ are connected,
and $x$ and $\rho(x)$ lie on opposite sides.
Therefore $\rho'$ interchanges these sides.
Consequently,
\[
\Bis(x,\rho(x))=\Fix\rho'\subset\Bis(u,v),
\]
which proves the reverse inclusion.
\qeds

The map $\rho$ as in the observation will be called a \emph{reflection} of $M$; the fixed-point set $\Bis(u,v)$ of $\rho$ will be called the \emph{mirror} of $\rho$ and denoted by $\Bis_\rho$.
If every mirror of $M$ has empty interior, then we say that $M$ has \emph{thin mirrors};
otherwise we say that it has \emph{thick mirrors}.

\begin{thm}{Lemma}\label{lem:rigidity}
Let $M$ be an all-set-homogeneous length space with thin mirrors.
Then
\begin{enumerate}[(i)]
\item no mirror of $M$ properly contains another;
\item an isometry of $M$ fixing a mirror pointwise and preserving its sides is the identity;
\item an isometry is determined by its values on two distinct mirrors;
\item the stabilizer of an oriented mirror acts transitively on that mirror.
\end{enumerate}

\end{thm}

\parit{Proof.}
For nested mirrors, their two sides and their common boundaries force equality.

Suppose $q=\iota(p)\ne p$ for an isometry $\iota$ of $M$ that fixes a mirror $\Bis_\rho$ pointwise.
Then $\Bis_\rho\subset \Bis(p,q)$ and therefore equality holds; in this case, $p$ and $q$ lie on the opposite sides of $\Bis_\rho$.

For transitivity on $\Bis_\rho$, choose points on one side of $\Bis_\rho$ sufficiently close to distinct points $x,y\in\Bis_\rho$, join them on this side, and choose a point $z$ with $|x-z|=|y-z|$.
The reflection swapping $x$ and $y$ fixes $z$ and $\rho(z)$, hence preserves their bisector and its chosen side.
\qeds

\begin{thm}{Lemma}\label{prop:thick mirror}
Any all-set-homogeneous length space with a thick mirror has local cut points.
\end{thm}

\parit{Proof.}
Let $\Bis_\rho$ be a thick mirror in an all-set-homogeneous length space $M$.
Choose a connected component $W\subset\Int_M(\Bis_\rho)$ and a path $\gamma\colon[0,1]\to M$
such that $\gamma(t)\in W$ for any $t<1$ and $\gamma(1)\in\partial W$.
Let $p=\gamma(0)$, $q=\gamma(1)$, and $s=|p-q|$.

Choose $z\in M\setminus \Bis_\rho$ sufficiently close to $q$,
and let $z'=\rho(z)$.
Observe that $|z-\gamma(t)|=|z'-\gamma(t)|$ for any $t$.
All-set-homogeneity implies that given $\eps>0$ there is $\delta>0$ such that if $|z-z'|<\delta$ and $\eps<|x-z|=|x-z'|<\tfrac s2$, then
\[\Ball(x,\delta)\subset \Bis(z,z')=\Bis_\rho;\]
the last equality follows from \ref{lem:reflection}.

It follows that $q$ is an isolated point of $\partial W$ and hence is a local cut point of~$M$.
\qeds

The following lemma is about isosceles-free spaces, which were studied by Christian Bargetz, Adam Bartoš, Wiesław Kubiś, and Franz Luggin \cite{bargetz-bartos-kubis-luggin}.

\begin{thm}{Lemma}\label{lem:equilateral-set}
Let $M$ be an all-set-homogeneous length space without local cut points.
Then any mirror $\Bis_\rho$ in $M$ contains a triple of distinct points $o$, $u$ and $v$ such that $|o-u|=|o-v|$.
\end{thm}

\parit{Proof.}
By \ref{prop:thick mirror}, $M$ has thin mirrors.
Furthermore, no mirror in $M$ has an isolated point; otherwise such a point would be a local cut point.

Suppose that the required triple of points in $\Bis_\rho$ does not exist.
Then an isometry of $\Bis_\rho$ would be determined by the image of any one point.

Fix $o\in \Bis_\rho$ and, for $x\in \Bis_\rho\setminus\{o\}$, let $\tau_x$ be the reflection swapping $o$ with~$x$.
The reflections $\tau_x$ preserve $\Bis_\rho$ and its sides, and their restrictions, together with the identity, form a group;
since $\tau_x^2=\id_M$ it must be commutative.
If $x\ne y$, then the composition $\tau_y\circ\tau_x$ must be another reflection, which is impossible: its fixed set is the intersection of their mirrors.
\qeds

\begin{thm}{Lemma}\label{lem:cont-on-mirror}
Let $M$ be an all-set-homogeneous length space with thin mirrors.
Consider two converging sequences of points $p_n\to p_\infty$, $q_n\to q_\infty$ in $M$ such that $p_\infty\ne q_\infty$, and $p_n\ne q_n$ for each $n$.
Suppose $\rho_\infty$ is the reflection swapping $p_\infty$ with $q_\infty$ and $\rho_n$ is the reflection swapping $p_n$ with $q_n$ for each $n$.
Then $\rho_n(z)\to z$ as $n\to\infty$ for any $z\in \Bis_{\rho_\infty}$.
\end{thm}

Since $z\in \Bis_{\rho_\infty}$, we have $\rho_\infty(z)=z$.
Once the main theorem is proved, we will also have that $\rho_n(z)\to \rho_\infty(z)$ for any $z$.

\parit{Proof.}
Fix $z\in\Bis_{\rho_\infty}$.
Since the mirror is thin, we can choose $x\in M\setminus\Bis_{\rho_\infty}$ arbitrarily close to $z$;
set $y=\rho_\infty(x)$.

The points $x$ and $y$ lie on opposite sides of $\Bis_{\rho_\infty}$;
that is, the values
\[
|x-p_\infty|-|x-q_\infty|
\quad\text{and}\quad
|y-p_\infty|-|y-q_\infty|
\]
have opposite signs.
Since $p_n\to p_\infty$ and $q_n\to q_\infty$,
the values
\[
|x-p_n|-|x-q_n|
\quad\text{and}\quad
|y-p_n|-|y-q_n|
\]
have opposite signs for all large $n$.
Therefore, an almost shortest path from $x$ to $y$ meets $\Bis_{\rho_n}$ at some point $w_n$.
Since $\rho_n(w_n)=w_n$, we have
\[
|z-\rho_n(z)|
\leqslant 2\cdot|z-w_n|,
\]
and the result follows.
\qeds

\section{Geodesic spaces}\label{Locally unique geodesics}

\subsection*{Splitting point}

Given two points $p$ and $q$ in a geodesic space $M$, let us denote by $\vv pq\ww$ the union of all geodesics from $p$ to $q$.
In other words,
\[x\in \vv pq\ww\quad\Leftrightarrow\quad|p-x|+|q-x|=|p-q|.\]

\begin{thm}{Lemma}\label{lem:split-time}
Let $M$ be an all-set-homogeneous geodesic space.
For any triple of points $p,q_1,q_2\in M$ there is a point $x\in \vv pq_1\ww\cap \vv pq_2\ww$ that maximizes the distance $|p-x|$.

\end{thm}

\parit{Proof.}
Let $r\df\sup\set{|p-x|}{x\in \vv pq_1\ww\cap \vv pq_2\ww}$.

Choose a sequence $x_1,x_2,\ldots\in \vv pq_1\ww\cap \vv pq_2\ww$ so that $r_n\df|p-x_n|$ is an increasing sequence converging to $r$.

Applying induction, we can assume
\[[px_1]\subset[px_2]\subset\dots
\eqlbl{one-star}\]
Indeed, suppose $[px_1]\subset\ldots\subset[px_n]$.
Choose a geodesic $[px_{n+1}]$ and a point $x'_n\in [px_{n+1}]$ such that $|p-x'_n|=|p-x_n|$.
Since $M$ is all-set-homogeneous, there is an isometry of $M$ that fixes points $p$, $q_1$, $q_2$ and sends $x'_n$ to $x_n$.
Replacing $x_{n+1}$ by its image under this isometry, we have $[px_n]\subset [px_{n+1}]$ for an appropriate choice of $[px_{n+1}]$.

By \ref{obs:geodesically-complete}, the nested union of geodesics in \ref{one-star} extends to a local geodesic defined on $\RR$, which contains the limit point $x$ of the sequence; hence the result.
\qeds

\subsection*{Splitting function}

Consider three distinct points $p$, $q_1$, and $q_2$ in an all-set-homogeneous geodesic space~$M$;
let $x$ be a splitting point provided by \ref{lem:split-time}.
If $|p-q_1|\z=|p-q_2|$ and $x\ne p$, then a triple of geodesics $[xq_1]$, $[xq_2]$, and $[xp]$ will be called a \emph{fork}.
Note that the unions $[xp]\cup [xq_1]$ and $[xp]\cup [xq_2]$ form geodesics $[pq_1]$ and $[pq_2]$.
Let $\gamma_0$, $\gamma_1$, and $\gamma_2$ be the arc-length parametrizations of $[xp]$, $[xq_1]$, and $[xq_2]$ starting from~$x$;
so, $x\z=\gamma_0(0)\z=\gamma_1(0)\z=\gamma_2(0)$.

\begin{wrapfigure}{o}{40mm}
\centering
\vskip-0mm
\includegraphics{mppics/pic-50}
\vskip-0mm
\end{wrapfigure}

Let us restrict each $\gamma_i$ to a segment $[0,\ell]$ for some $0<\ell<\tfrac12\cdot\diam M$.
Note that $\gamma_1(t)=\gamma_2(t)$ only if $t=0$.
In particular, \[s(t)\df|\gamma_1(t)-\gamma_2(t)|>0
\quad\text{if}\quad
t>0.\]
The function $s$ will be called the \emph{splitting function} of the fork.

\begin{thm}{Lemma}\label{lem:splitting-function}
Let $s$ be a splitting function of a fork in an all-set-homogeneous geodesic space $M$.
Then $s$ is strictly increasing.
Moreover, the splitting function is essentially independent of the choice of fork;
that is, if $s$ and $s'$ are splitting functions for different forks in $M$,
then $s$ coincides with $s'$ in the common interval of definition.

\end{thm}

\parit{Proof.}
Clearly $s(0)=0$;
by the triangle inequality, $s$ is $2$-Lipschitz.

Observe that $s$ is strictly increasing.
Indeed, suppose $s(r_0)\geqslant s(r_1)$ for some $r_0<r_1$.
By the intermediate value theorem, we can find $0<r\leqslant r_0<r'\leqslant r_1$
 such that $s(r)=s(r')$.
If $\gamma_0$, $\gamma_1$, and $\gamma_2$ are as above,
then
\[\gamma_1(r)\mapsto\gamma_1(r'),
\quad
\gamma_2(r)\mapsto\gamma_2(r'),
\quad
\gamma_0(r')\mapsto\gamma_0(r)\]
defines a distance-preserving map on a 3-point set,
which can be extended to an isometry of $M$.
The maximal splitting distances for the two triples are $r'$ and $r$, and the isometry must preserve this quantity.
Therefore $r'=r$ --- a contradiction.
Here we apply \ref{lem:split-time} to two triples:
$\gamma_0(r')$, $\gamma_1(r)$, $\gamma_2(r)$
and $\gamma_0(r)$, $\gamma_1(r')$, $\gamma_2(r')$.

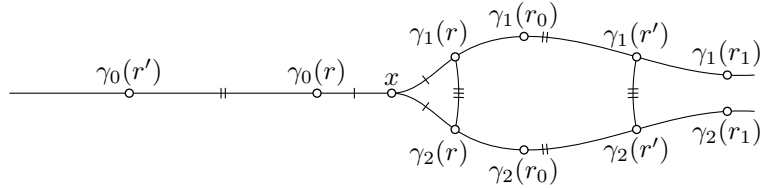
\begin{figure}[h!]
\vskip-0mm
\centering
\includegraphics{mppics/pic-20}
\end{figure}

The same argument proves the second statement.
Indeed, assume $s'$ is a splitting function for another fork described by geodesics $\gamma'_0$, $\gamma'_1$, and $\gamma'_2$.
If $s(r)=s'(r')$ for $r$ and $r'$ in the common domain, then
\[\gamma'_1(r')\mapsto\gamma_1(r),\quad
\gamma'_2(r')\mapsto\gamma_2(r),\quad
\gamma'_0(r)\z\mapsto\gamma_0(r')\]
defines a distance-preserving map, which extends to an isometry of $M$, implying $r=r'$.
Since both functions are continuous and strictly increasing from $0$, it follows that $s=s'$ on their common domain.
\qeds

The lemma together with homogeneity implies the following.

\begin{thm}{Corollary}\label{cor:splitting-function}
Let $M$ be an all-set-homogeneous geodesic space, and let $s$ be its splitting function.
Suppose $\gamma_1,\gamma_2\colon[0,r]\to M$ are two geodesics such that $\gamma_1(0)=\gamma_2(0)$ and $|\gamma_1(r)-\gamma_2(r)|=s(r)$.
Then
\[|\gamma_1(r')-\gamma_2(r')|= s(r')
\quad\text{for any}\quad
r'\in[0,r].\]

\end{thm}

\subsection*{Rhombus}

\begin{thm}{Lemma}\label{lem:rhombus}
Let $M$ be a geodesic all-set-homogeneous space with locally non\-unique geodesics.
Then for any sufficiently small $r$, there is a configuration of four points $p_1$, $p_2$, $q_1$, and $q_2$ such that
\[|p_1-p_2|=2\cdot r,\quad |q_1-q_2|=s(r),\quad\text{and}\quad |p_i-q_j|=r\]
for all $i$ and $j$;
here $s$ denotes the splitting function of $M$.
\end{thm}

Such a configuration of points will be called an \emph{$r$-rhombus}.

\parit{Proof.}
Suppose $0<R<\diam M$.

Since the geodesics are not locally unique, there is a configuration of distinct points $v_1,v_2,q_1,q_2$ such that $|v_1-v_2|$ is small, $q_1,q_2\in \vv v_1v_2\ww$ and $|v_1-q_1|\z=|v_1-q_2|$ (and hence $|v_2-q_1|=|v_2-q_2|$).
Since the space is geodesic and all-set-homogeneous, points $v_1$ and $v_2$ lie on an $R$-long geodesic $[w_1w_2]$,
both close to its midpoint.

Note that $|q_1-q_2|\leqslant|v_1-v_2|$;
hence it is small.
Therefore, $|q_1-q_2|=s(r)$ for some $r\ll R$.
Let $p_i$ be the point provided by \ref{lem:split-time} for $w_i$, $q_1$ and $q_2$.
Note that $p_1$, $p_2$, $q_1$, $q_2$ form an $r$-rhombus.

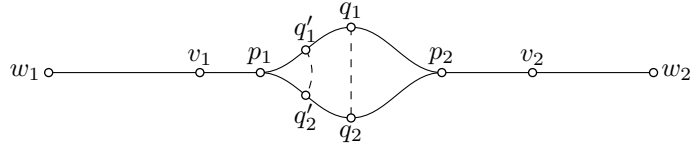
\begin{figure}[b!]
\vskip-0mm
\centering
\includegraphics{mppics/pic-30}
\end{figure}

Choose a pair of points $q'_1\in [p_1q_1]$ and $q'_2\in [p_1q_2]$ such that $|p_1-q'_1|\z=|p_1-q'_2|=r'\leqslant r$.
By \ref{cor:splitting-function}, $|q'_1-q'_2|=s(r')$.
Repeating the construction above, we get an $r'$-rhombus for any $r'\in(0,r]$.
\qeds

\subsection*{Self-improvement}

\begin{thm}{Lemma}
Let $M$ be an all-set-homogeneous geodesic space with locally nonunique geodesics.
Suppose $r>0$ is sufficiently small.
Then there are $p,q\z\in M$ and an infinite sequence
$x_1,x_2,\ldots\in M$ such that
\[|p-q|=2\cdot r,
\quad
|p-x_i|=|q-x_i|= r,
\quad\text{and}\quad
|x_i-x_j|=s(r)\] for $ i\ne j$.
\end{thm}

\begin{figure}[h!]
\centering\vskip-0mm
\includegraphics{mppics/pic-40}
\vskip-0mm
\end{figure}

\parit{Proof.}
In the proof we construct and use the configuration of points shown in the picture.
Each marked segment has length $r$ times its Roman numeral, and the black points are pairwise at distance $s(r)$.

Choose two points $p$ and $q$ such that $|p-q|=2\cdot r$.
Since $r$ is small, we can extend this pair to a $2\cdot r$-rhombus $p,p',q,q'\in M$;
that is, $|p-p'|=4\cdot r$, $|q-q'|=s(2\cdot r)$, and $|p-q|=|p'-q|=|p-q'|=|p'-q'|=2\cdot r$.
Furthermore, there is a point $z$ such that $p\in \vv zp'\ww$ and $|z-p|=2\cdot r$.

Suppose we have a finite sequence $x_1,\ldots,x_n$ that meets the conditions.
If we construct one more point $x_{n+1}$, then induction finishes the proof.

Let $y$ be the midpoint of $[pq']$.
By Corollary~\ref{cor:splitting-function}, $|y-x_i|=s(r)$ for any~$i$.
By the triangle inequality, $|p'-y|=|z-y|=|p'-x_i|=|z-x_i|=3\cdot r$.
Thus an isometry $\iota$ of $M$ maps $z$ to $p'$ and fixes $y,x_1,\ldots,x_n$.
Set $w=\iota(p)$;
the triangle inequality gives
\[|p-w|=2\cdot r
\quad\text{and}\quad
|p-y|=|w-y|=|p-x_i|=|w-x_i|= r\]
for any $i$.
Therefore, there is an isometry $\iota'$ that fixes $p,x_1,\ldots,x_n$ and moves $w$ to $q$.
Note that the point $x_{n+1}=\iota'(y)$ meets all the requirements.
\qeds

\begin{thm}{Proposition}\label{prop:unique-geod}
All-set-homogeneous geodesic spaces have locally unique geodesics.

\end{thm}

\parit{Proof.}
Assume, to the contrary, that $M$ is an all-set-homogeneous geodesic space with locally nonunique geodesics.
Choose two sufficiently close points $p,q\z\in M$;
let $r=\tfrac12\cdot|p-q|$.
By the lemma, the set of midpoints of $p$ and $q$ contains an infinite $s(r)$-equilateral set,
which contradicts \ref{obs:spheres} and~\ref{obs:equilateral}.
\qeds

\subsection*{Universal tree}

\begin{thm}{Proposition}\label{prop:loc-tree}
Let $M$ be an all-set-homogeneous geodesic space with locally unique geodesics.
If it has bifurcating geodesics, then $M$ is a universal metric tree of finite valence.
\end{thm}

\parit{Proof.}
Homogeneity and local uniqueness of geodesics imply that there is $R>0$ such that if $|x-y|<R$, then $[xy]$ is uniquely defined.

Let $s$ be the splitting function of $M$; it is defined since $M$ has bifurcating geodesics.%
\footnote{For metric trees we must have $s(t)\equiv2\cdot t$.}

Fix a point $p\in M$ and $r>0$.
Given a point $x$, denote by $x_r$ the point on $[px]$ such that $|p-x_r|=r$;
it is uniquely defined if $r\leqslant |p-x|<R$.
Observe that for any small $\eps>0$ there is $\delta>0$ such that
\[|p-x|,|p-y|\in (r+\eps,R)\  \text{and}\ |x-y|\in(0,\delta)\ \Rightarrow\  x_r=y_r.\]

\begin{wrapfigure}{o}{40mm}
\centering
\vskip-0mm
\includegraphics{mppics/pic-60}
\vskip-0mm
\end{wrapfigure}

Indeed, let $a=|p-x|$;
we can assume that $a\leqslant |p-y|$.
By \ref{lem:splitting-function}, we can choose $\delta>0$ such that $s(\eps)>2\cdot \delta$.
By the triangle inequality,
\[|x-y|<\delta\quad\Rightarrow\quad|x-y_a|<2\cdot\delta.\]
Suppose $w$ is the splitting point of $[px]$ and $[py]$, and let $b=|w-x|=|w-y_a|$.
Then, by 3-point-homogeneity, $s(b)\leqslant|x-y_a|$,
and hence $b<\eps$ and $x_r=y_r$;
here we assume that $s$ is defined on $[0,R]$.

Choose a geodesic triangle $[pxy]$ in $M$ such that $|p-x|<\tfrac R2$, $|p-y|<\tfrac R2$;
let $r=\tfrac12\cdot(|p-x|+|p-y|-|x-y|)$.
Let us show that $[pxy]$ is a tripod.

If not, then $|p-z|>r$ for any $z\in [xy]$;
note also that $|p-z|<R$.
From above, $z\mapsto z_r$ is locally constant and therefore constant on~$[xy]$.
It follows that $x_r=y_r\in [px]\cap[py]$.
Since $|x-x_r|+|y-y_r|=|x-y|$, we have $[xy]\subset [px]\cup[py]$ --- a contradiction.

Thus $M$ meets the assumptions in \ref{prop:global-tree}, which finishes the proof.
\qeds

Combining the last result with Birkhoff's theorem, we get the following.

\begin{thm}{Corollary}\label{thm:tree}
Any all-set-homogeneous geodesic space with locally unique geodesics is a universal tree of finite valence or a classical space.
\end{thm}

\section{Length spaces}

\subsection*{Spherical curve}

\begin{thm}{Monster lemma}\label{clm:rectifiable-curve}
Let $M$ be an all-set-homogeneous length space with thin mirrors.
Suppose there are two distinct commuting reflections $\rho$ and $\sigma$ of $M$ such that
\[A\df\Bis_\rho\cap\Bis_\sigma\ne\emptyset.\]
Then there are arbitrarily small radii $r>0$ such that a sphere $\Sph(o,r)$ in $M$ contains a nontrivial rectifiable curve.
\end{thm}

The proof is elementary, but delicate: one false step sends us over the edge.
The following three statements will follow from the main theorem, but only after it is proved.
Therefore they cannot be assumed; keep this in mind while reading further.
\begin{enumerate}[(i)]
\item If two distinct mirrors contain $A$, then they have no other points of intersection.
\item If two distinct mirrors in $M$ intersect, then one mirror intersects both open sides of the other.
\item Let $\rho_1,\rho_2,\dots$, and $\rho_\infty$ be reflections of $M$.
Suppose $\rho_n$ swaps distinct points $p_n$ with $q_n$ for each $n$;
furthermore, $p_n\to p_\infty$, $q_n\to q_\infty$, $p_\infty\ne q_\infty$, and
$\rho_\infty$ swaps $p_\infty$ with $q_\infty$.
Then $\rho_n(z)\to \rho_\infty(z)$ for any $z\in M$.
\end{enumerate}
If the last statement is known, then the monster lemma can be proved in a few lines; instead we have to use its weaker version \ref{lem:cont-on-mirror}.

\parit{Proof.}
Let us denote by $\bar M$ the completion of $M$; any isometry of $M$ will be identified with its extension to $\bar M$.

Consider a nested sequence of dihedral groups $D_2<D_4<D_{8}<\ldots$;
denote by $D_{\mathrm{dyad}}$ its union.
We construct an action of $D_{\mathrm{dyad}}$ on $M$ that mimics the action of a dense subgroup of $\OO(2)$, fixing the axis $A$ pointwise.
Further, we construct a rectifiable closed curve $\beta_\infty$ in $\bar M$ that is invariant with respect to the $D_{\mathrm{dyad}}$-action.
Finally, we show that $\beta_\infty$ lies in $M$; thus it gives the required rectifiable curve.
(The construction also leads to an $\OO(2)$-action on $M$ fixing $A$ pointwise, with $\beta_\infty$ as an orbit, but we will not use this.)

\parit{Path between mirrors.}
By \ref{lem:bisector}, $M\setminus A$ is path-connected.
Indeed, both path-connected components of $M\setminus \Bis_\rho$ intersect both path-connected components of $M\setminus \Bis_\sigma$.
Thus we can choose a rectifiable path $\alpha_1$ in $M\setminus A$ from $\Bis_\rho$ to $\Bis_{\sigma}$.

Moreover, given $o\in A$ we can assume that $\alpha_1$ lies arbitrarily close to $o$.
Indeed, choose a point $x\notin \Bis_\rho\cup \Bis_{\sigma}$ sufficiently close to $o$.
Connect $\rho(x)$ to $x$ and $x$ to $\sigma(x)$ by almost shortest paths;
note that they cross $\Bis_\rho$ and $\Bis_{\sigma}$ and do not intersect $\Bis_{\sigma}$ and $\Bis_\rho$, respectively.
Thus the concatenation of these two paths lies in $M\setminus A$ and contains points on $\Bis_\rho$ and $\Bis_{\sigma}$.
By the triangle inequality the path lies in a small neighborhood of $o$.

\parit{Dyadic action.}
Let us recursively construct a sequence of reflections $\sigma=\tau_1$, $\tau_2,\ldots$ and paths $\alpha_1$, $\alpha_2,\ldots$ in $M\setminus A$
such that $\alpha_n$ runs from $\Bis_\rho$ to $\Bis_{\tau_n}$,
\[(\tau_{n+1}\circ\rho)^2=\tau_n\circ\rho,
\quad
\Fix (\tau_n\circ\rho)=A,
\quad\text{and}\quad
2\cdot \length \alpha_{n+1}\leqslant \length \alpha_{n}
\eqlbl{eq:sigma-n-rho}\]
for each $n$.

\begin{figure}[ht!]
\centering\vskip-0mm
\includegraphics{mppics/pic-90}
\vskip-0mm
\end{figure}

Suppose $\tau_n$ and $\alpha_n$ are constructed.
By the intermediate value theorem there is a point $x\in \alpha_n$ such that $|x-\rho(x)|=|x-\tau_n(x)|$.
Let $\tau_{n+1}$ be the reflection swapping $\rho(x)$ and $\tau_n(x)$; note that $\tau_{n+1}$ fixes $x$ and each point in $A$.
Now $\tau_{n+1}\circ\rho\circ\tau_{n+1}=\tau_n$ since both reflections swap $x$ with $\tau_n(x)$, and we have the first identity in \ref{eq:sigma-n-rho}.

Put $\theta_n=\tau_n\circ\rho$.
Note that $\Fix \theta_1=\Bis_{\tau_1}\cap \Bis_\rho=A$.
Suppose
\[\Fix \theta_n=\Bis_{\tau_n}\cap \Bis_\rho=A.
\eqlbl{eq:Fix=A}\]
By the first identity in \ref{eq:sigma-n-rho}, $\theta_{n+1}^2=\theta_n$, and we have $\Fix \theta_{n+1}\subset \Fix \theta_n$.
Since $\tau_{n+1}$ fixes each point in $A$, we have the opposite inclusion.
Since
\[\Fix \theta_{n+1}\supset \Bis_{\tau_{n+1}}\cap \Bis_\rho\supset A,\]
we get \ref{eq:Fix=A} for $n+1$.

Now choose a shortest subarc of $\alpha_n$ that runs from $\Bis_\rho$ to $\Bis_{\tau_{n+1}}$ or from $\Bis_{\tau_{n+1}}$ to $\Bis_{\tau_n}$;
reflect it by $\tau_{n+1}$ if necessary, so it runs from $\Bis_\rho$ to $\Bis_{\tau_{n+1}}$ and denote the resulting path by $\alpha_{n+1}$.

The isometries $\rho$ and $\tau_n$ generate the dihedral  group $D_{2^{n}}$.
It mimics reflections and rotations with the axis $A$;
the isometry $\theta_n$ corresponds to a rotation by angle $\pi/2^{n-1}$;
so it has order $2^n$.
The isometries $\rho$ and $\tau_n$ for all $n$ generate $D_{\mathrm{dyad}}$, which is the union of the nested sequence of subgroups $D_2<D_4<\dots$;
the group $D_{\mathrm{dyad}}$ is naturally isomorphic to a subgroup of $\OO(2)$.
If $\theta_n^k\ne \id_M$, then one of its powers is $\theta_1$; therefore, $\theta_n^k$ fixes precisely $A$.
Furthermore the mirrors of any two distinct reflections in $D_{\mathrm{dyad}}$ intersect precisely at~$A$.

\parit{Spherical curve in $\bar M$.}
The arc $\alpha_n$ has $2^{n+1}$ copies obtained by applying $\rho$, $\tau_n$, and their compositions.
The resulting closed curve $\beta_n$ has $D_{2^{n}}$-symmetry;
furthermore, the lengths of all $\beta_n$ are uniformly bounded.
We will assume that $\beta_n$ comes with a $D_{2^{n}}$-symmetric parametrization by $\S^1=\RR/\ZZ$ starting at $\Bis_\rho$ with Lipschitz constant $4\,\cdot\length\alpha_1$.

Note that for compatible parametrizations, we have
\[|\beta_{n+1}(t)-\beta_n(t)|\leqslant \mathrm{const}/2^n.\]
Hence $\beta_n$ converge to a curve $\beta_\infty$ in $\bar M$;
the limit parametrization respects the $D_{\mathrm{dyad}}$-action.

\parit{Limit action.}
Denote by $\tilde\alpha_n$ the subarc of $\alpha_1$ that corresponds to $\alpha_n$ after unfolding it back to $\alpha_1$.
Let $p_n$ be the endpoint of $\tilde\alpha_n$ that corresponds to the endpoint of $\alpha_n$ on $\Bis_\rho$.
Then $p_n\to p_\infty\in\alpha_1\subset M\setminus A$.
Note that $p_\infty$ lies in $\beta_\infty$ and hence the whole $D_{\mathrm{dyad}}$-orbit of $p_\infty$ lies on $\beta_\infty$.

Let $\rho_n$ be the corresponding conjugate of $\rho$; so, $\rho_n(p_n)=p_n$.
It has a partner reflection $\sigma_n=\theta_1\circ \rho_n$, which swaps $p_n$ with $\theta_1(p_n)$.
Let $\sigma_\infty$ be the reflection swapping $p_\infty$ and $\theta_1(p_\infty)$.
By \ref{lem:cont-on-mirror}, we have
\[
z\in\Bis_{\sigma_\infty}
\quad\Rightarrow\quad
\sigma_n(z)\to z
\quad\Rightarrow\quad
\rho_n(z)\to \theta_1(z).
\]

A path from $p_\infty$ to $\theta_1(p_\infty)$ in $M\setminus A$
gives a point $q_\infty\in\Bis_{\sigma_\infty}\setminus A$.
By \ref{lem:cont-on-mirror} and \ref{lem:reflection}, we can choose $q_n\in\Bis_{\sigma_n}$ such that $q_n\to q_\infty$.
Note that $\rho_n(q_n)=\theta_1(q_n)$ for all large $n$.


Let $\rho_\infty$ be the reflection swapping $q_\infty$ with $\theta_1(q_\infty)$.
By \ref{lem:cont-on-mirror}, $\rho_n(z)\to z$ for every $z\in\Bis_{\rho_\infty}$.
Together with the convergence on $\Bis_{\sigma_\infty}$, this implies
\[
|z-w|
=\lim_{n\to\infty}|\rho_n(z)-\rho_n(w)|
=|z-\theta_1(w)|
\]
for $z\in\Bis_{\rho_\infty}$ and $w\in\Bis_{\sigma_\infty}$.
Consequently, if $w\in\Bis_{\sigma_\infty}\setminus A$, then $\Bis_{\rho_\infty}\z\subset\Bis(w,\theta_1(w))$.
By \ref{lem:rigidity}, these mirrors coincide and hence
\[\rho_\infty=\theta_1\circ\sigma_\infty;\]
in particular, $\theta_1\circ\sigma_\infty$ is a reflection.

Let us show that
\[
\rho_\infty\circ\theta_m\circ\rho_\infty=\theta_m^{-1}
\eqlbl{eq:rhothetarhotheta}
\]
for every $m$.

Fix $m$;
for $n\geqslant m$, we have $\rho_n\circ\theta_m=\theta_m^{-1}\circ\rho_n$.
From above
\[
\rho_n(z)\to\rho_\infty(z),
\quad\text{and therefore}\quad
\rho_n\circ\theta_m(z)
\to
\theta_m^{-1}\circ\rho_\infty(z),
\]
for every $z\in\Bis_{\rho_\infty}\cup\Bis_{\sigma_\infty}$;
indeed, if $z\in\Bis_{\rho_\infty}$, then $\rho_\infty(z)=z$,
and if $z\in\Bis_{\sigma_\infty}$, then $\rho_\infty(z)=\theta_1(z)$.
Thus the maps $\rho_n$ converge on $\Bis_{\rho_\infty}\cup\Bis_{\sigma_\infty}$ and its image under $\theta_m$ to a distance-preserving map with values in~$M$.
By all-set-homogeneity, this partial isometry extends
to an isometry $\iota$ of $M$, which agrees with $\rho_\infty$ on $\Bis_{\rho_\infty}\cup\Bis_{\sigma_\infty}$;
by \ref{lem:rigidity}, $\iota=\rho_\infty$.
Therefore $\rho_\infty\circ\theta_m$ and
$\theta_m^{-1}\circ\rho_\infty$ agree on $\Bis_{\rho_\infty}\cup\Bis_{\sigma_\infty}$.
Another application of \ref{lem:rigidity} provides \ref{eq:rhothetarhotheta}.

For every $m\geqslant1$ and $k\in\ZZ$, we have
\[
\theta_m^k\circ\rho_\infty
=
\theta_{m+1}^k\circ\rho_\infty\circ\theta_{m+1}^{-k},
\]
so $\theta_m^k\circ\rho_\infty$ is conjugate to $\rho_\infty$;
in particular, $\theta_m^k\circ\rho_\infty$ is a reflection.
The partner reflection $\sigma_\infty=\theta_1\circ\rho_\infty$ belongs to the same family.

Let us reparametrize $\beta_\infty$ starting from $p_\infty$;
so
\[
\beta_\infty(0)=p_\infty,
\quad
\beta_\infty(s+\tfrac1{2^{m}})=\theta_m\circ\beta_\infty(s),
\quad\text{and}\quad
\beta_\infty(-s)=\rho_\infty\circ\beta_\infty(s).
\]
From now on $\beta_\infty(d)\in M$ for any dyadic $d\in\RR/\ZZ$.

\parit{Pushing $\beta_\infty$ into $M$.}
Let us show that $\beta_\infty(t_\infty)\in M$ for any $t_\infty\in\RR/\ZZ$.

For every dyadic $t=\tfrac k{2^n}$, let $\Theta_t=\theta_n^k$.
Choose dyadic numbers $t_j\to t_\infty$.
Let $\eta_j=\Theta_{t_j}\circ\rho_\infty$; so $\eta_j\circ\beta_\infty(s)=\beta_\infty(t_j-s)$.

Choose a rectifiable path $\gamma$ from $p_\infty$ to $\theta_1(p_\infty)$ in $M\setminus A$.
Since $\theta_1$ interchanges the sides of every $\Bis_{\eta_j}$, we can choose $y_j\in\gamma\cap\Bis_{\eta_j}$.
After passing to a subsequence, $y_j\to y_\infty\in\gamma\subset M\setminus A$.
Let us apply the preceding construction of a limiting pair of reflections to $\eta_j$ and $\theta_1\circ\eta_j$,
where $y_\infty$ plays the role of $p_\infty$.
It gives a reflection $\eta_\infty$ such that $\eta_j(z)\to\eta_\infty(z)$
for every $z\in\Bis_{\eta_\infty}$, and $\eta_\infty\circ\theta_m=\theta_m^{-1}\circ\eta_\infty$ for every $m$.
Moreover, $\theta_1\circ\eta_\infty$ is a reflection.

The mirror $\Bis_{\eta_\infty}$ meets every closed curve $\beta_n$:
indeed, $\beta_n$ is $\theta_1$-invariant,
and $\theta_1$ interchanges the sides of $\Bis_{\eta_\infty}$.
Choose $z_n\in\beta_n\cap\Bis_{\eta_\infty}$.
Since $\beta_n$ converge to the closed curve $\beta_\infty$, a subsequence of $z_n$ converges to a point $z_\infty$ in $\beta_\infty\cap\Cl_{\bar M}\Bis_{\eta_\infty}$.
Thus $\eta_\infty(z_\infty)=z_\infty$; recall that each isometry of $M$ is identified with its extension to $\bar M$.

By \ref{lem:cont-on-mirror}, we have convergence $\eta_j\to\eta_\infty$ on $\Bis_{\eta_\infty}$,
and it also holds at $z_\infty$.
Indeed, for $z\in\Bis_{\eta_\infty}$,
\[
\limsup_{j\to\infty}|\eta_j(z_\infty)-\eta_\infty(z_\infty)|
\leqslant 2\cdot|z_\infty-z|,
\]
and the right-hand side can be made arbitrarily small.

Choose $s$ such that $z_\infty=\beta_\infty(s)$.
Since $\eta_j(z_\infty)\to z_\infty$,
we have $\beta_\infty(t_\infty-s)\z=z_\infty$.
If $d$ is dyadic, then
\begin{align*}
\eta_\infty\circ\beta_\infty(s+d)
&=\eta_\infty\circ\Theta_d(z_\infty)=
\\
&=\Theta_{-d}\circ\eta_\infty(z_\infty)=
\\
&=\Theta_{-d}(z_\infty)=
\\
&=\beta_\infty(t_\infty-s-d).
\end{align*}
By continuity, the same holds for every $d$.
If $s+d=0$, we obtain
\[
\beta_\infty(t_\infty)=\eta_\infty\circ\beta_\infty(0)=\eta_\infty(p_\infty)\in M.
\]
Since $t_\infty$ was arbitrary, $\beta_\infty$ lies in $M$.

Thus $\beta_\infty$ is a rectifiable curve in $M$,
which is not constant since $p_\infty\ne \theta_1(p_\infty)$.
Since the $D_{\mathrm{dyad}}$-action fixes $A$ pointwise and has a dense orbit in $\beta_\infty$, the curve lies in a sphere centered at any point $o\in A$.
Furthermore, the radius of the sphere is $r=|o-p_\infty|\z=|o-\alpha_1(t)|$ for some $t$;
in particular, $r$ can be taken arbitrarily small.
\qeds

\subsection*{Completeness}\label{sec:completeness}

\begin{thm}{Proposition}\label{prop:complete}
Any all-set-homogeneous length space is complete.
\end{thm}

\parit{Proof of \ref{prop:complete}.}
Arguing by contradiction, assume $M$ is a noncomplete all-set-homogeneous length space.
By \ref{prop:cut-point} and \ref{prop:thick mirror}, $M$ has thin mirrors.
Furthermore, by \ref{lem:equilateral-set}, we can find two reflections of $M$, say $\rho$ and $\sigma$, that meet the assumptions of \ref{clm:rectifiable-curve}; the reflection $\sigma$ swaps $u$ with $v$ provided by \ref{lem:equilateral-set}.

Choose a point $\xi\in \bar M\setminus M$; recall that $\bar M$ denotes the completion of $M$.

Choose a sufficiently small radius $r_1>0$ supplied by \ref{clm:rectifiable-curve} and a point $p_1\in M$ such that $|p_1-\xi|_{\bar M}=r_1$;
let $\Sigma_1=\Sph(p_1,r_1)_M$.
Note that $\Sigma_1$ is not complete and it contains $\xi$ in its completion;
more precisely, $\Cl_{\bar M}\Sigma_1\ni \xi$.
Indeed, choose a path, say $\gamma$, that starts at $p_1$ and crosses $\Sigma_1$.
Then given a sequence of points $x_n\to \xi$, there is a sequence of points $y_n\in\gamma$ such that $|p_1-x_n|=|p_1-y_n|$ for any $n$ and $y_n\to y_\infty\in\gamma$ after passing to a subsequence.
By homogeneity, there is an isometry $\iota_n$ that fixes $p_1$ and moves $y_n$ to $x_n$, so $\iota_n(y_\infty)\in \Sigma_1$ for every $n$ and $\iota_n(y_\infty)\to \xi$.

Furthermore, by \ref{lem:length-sphere}, there is a connected component $\Upsilon_1$ of $\Sigma_1$ that contains $\xi$ in its completion.
Let us equip $\Upsilon_1$ with the induced length metric from $M$.
By \ref{lem:length-sphere} we can repeat the procedure in $\Upsilon_1$,
obtaining a sphere
\[\Sigma_2=\Sph(p_2,R_2)_{\Upsilon_1}=\Sph(p_1,r_1)_M\cap\Sph(p_2,r_2)_M;\]
here $R_2=\psi_{\Sigma_1}(r_2)$; we can assume that $r_2$ is much smaller than $r_1$ and $\eps_1$, where $\eps_1$ is provided by \ref{lem:length-sphere} for $\Sph(p_1,r_1)_M$.

Continuing this way, we get a sequence of spheres $\Sph(p_n,r_n)_M$ such that
\[\Sph(p_1,r_1)_M\cap\ldots\cap\Sph(p_{n+1},r_{n+1})_M\supset \Upsilon_{n+1}\ne \emptyset\]
for each $n$.
By \ref{obs:FIP},
$\Sph(p_1,r_1)_M\cap  \Sph(p_2,r_2)_M\cap\ldots\ne\emptyset$,
but since $|p_n-\xi|_{\bar M}=r_n$ and $r_n\to 0$ as $n\to\infty$, this intersection has to be empty --- a contradiction.
\qeds

\subsection*{Separated sequence}

Let $M$ be an all-set-homogeneous metric space.
Fix two positive real numbers $a$, $b$ and let $c=a+b$.
Consider the set $\Lambda_{a,c}$ of all $\lambda\geqslant 0$
such that $a=|u-v|$, $c=|u-w|$, and $b+\lambda=|v-w|$ for some points $u,v,w\in M$;
further assume $\Lambda_{a,c}\ne\varnothing$.

\begin{figure}[h!]
\centering\vskip-0mm
\includegraphics{mppics/pic-70}
\vskip-0mm
\end{figure}

For a pair $\lambda,\mu\in\Lambda_{a,c}$ we may consider the distance $|\lambda-\mu|$ induced from $\RR$,
but there are also two other distances $\theta(\lambda,\mu)$ and $\Theta(\lambda,\mu)$ defined as the greatest lower bounds of the values $\ell$ and $L$ for configurations of points shown in the picture.
Formally, they are defined as
\begin{align*}
\theta(\lambda,\mu)
&\df\inf\set{|v_1-v_2|}{
\substack{
|u-v_i|=a,\ |v_1-w|=b+\lambda,
\\
|u-w|=c,\ |v_2-w|=b+\mu.
}
},
\\
\Theta(\lambda,\mu)
&\df\inf\set{|w_1-w_2|}{
\substack{
|u-v|=a,\ |v-w_1|=b+\lambda,
\\
|u-w_i|=c,\ |v-w_2|=b+\mu.
}
},
\end{align*}
where $u,v,v_1,v_2,w,w_1,w_2\in M$.

The triangle inequalities for $\theta$ and $\Theta$ follow by applying $3$-point-homogeneity to identify the common triangle in two almost minimizing configurations. Evidently, $|\lambda-\mu|\leqslant \theta(\lambda,\mu)$ and $|\lambda-\mu|\leqslant \Theta(\lambda,\mu)$; hence $\theta$ and $\Theta$ are metrics on $\Lambda_{a,c}$.

If $M$ is complete, then so are $(\Lambda_{a,c},\theta)$ and $(\Lambda_{a,c},\Theta)$.
Indeed, pass to a subsequence of a Cauchy sequence with summable
successive distances.
By $3$-point-homogeneity, its terms can be represented by points
on a fixed sphere whose successive distances are summable.
Their limit represents the limit in the corresponding metric.

\begin{thm}{Lemma}\label{lem:cauchy-radius}
Let $M$ be a complete all-set-homogeneous length space.
Assume there are no points $u,v,w\in M$ such that
\[
|u-v|=a,\quad |v-w|=b,\quad |u-w|=c,
\]
where $a>0$, $b>0$, and $c=a+b<\diam M$.
Then any sequence in $\Lambda_{a,c}$ converging to $0$ contains a subsequence $\lambda_1,\lambda_2,\ldots$
such that
\[\theta(\lambda_i,\lambda_j)>\delta
\quad\text{and}\quad
\Theta(\lambda_i,\lambda_j)>\delta
\]
for some fixed $\delta>0$ and all $i\ne j$.
\end{thm}

\parit{Proof.}
Choose $u,w\in M$ with $|u-w|=c$ and a sequence $v_1,v_2,\ldots\in\Sph(u,a)$ such that $\lambda_n=|v_n-w|- b$ is the given sequence.

By assumption, $0\notin\Lambda_{a,c}$.
Therefore, the sequence $\lambda_1,\lambda_2,\ldots$ has no $\theta$-Cauchy subsequence.
Consequently, $\lambda_1,\lambda_2,\ldots$ has an infinite $\delta$-separated
subsequence for some $\delta>0$. The same argument applies to $\Theta$.
\qeds

\subsection*{Descent}

\begin{thm}{Lemma}\label{lem:plus-eps}
Let $M$ be an all-set-homogeneous complete length space.
Suppose $|u-v|<\diam M$.
Then there is a point $w$ such that
$|u-v|+|v-w|=|u-w|$ and $|v-w|>0$ is arbitrarily small.
\end{thm}

\parit{Proof.}
Let $a=|u-v|$; choose $r>0$ so that $a+r<\diam M$.
We can assume $a>0$; the case $a=0$ is immediate.
Suppose $b\in(0,r]$, and $c=a+b$.

Choose $w\in M$ such that $|u-w|=c$.
Let
\[h(b)=\lim_{\lambda\to 0+}\diam(\Sph(u,a)\cap \Ball(w,b+\lambda)).\]
By 2-point-homogeneity, $h(b)$ does not depend on the choice of $u$ and $w$.
Note that $h$ is nondecreasing and upper semicontinuous on $(0,r]$; furthermore,
\[
h(b)\leqslant2\cdot b.
\eqlbl{eq:sigma<2b}
\]

Suppose there is no point $w$ with $0<|v-w|\leqslant r$
that meets the requirements;
in other words, $M$ has no triangle with side lengths
$a$, $b$, and $c=a+b$  for $b\in(0,r]$.
By \ref{lem:cauchy-radius}, $h(b)>0$ for any $b\in(0,r]$.

We will show that \textit{$h$ is constant on $(0,r]$} and arrive at a contradiction with~\ref{eq:sigma<2b}.
Note that it is sufficient to show that for each $b\in(0,r]$
there is $b'\in(0,b)$ such that
\[
h(b')\geqslant h(b).
\eqlbl{sigma=const}
\]

Given $b\in(0,r]$ there is $\delta>0$ such that for any $\eps>0$ there are points $x,y,z\in \Sph(u,a)$ such that
\[\begin{aligned}
\eps&> h(b)-|y-z|,
&
\lambda(x)&>\lambda(y),
&
\Theta(\lambda(x),\lambda(y))&>\tfrac\delta2,
\\
\eps&>\lambda(x),
&
\lambda(x)&>\lambda(z),
&
\Theta(\lambda(x),\lambda(z))&>\tfrac\delta2,
\end{aligned}
\eqlbl{eq:lambda-eps}\]
where $\lambda(x)\df|w-x|-b$.

Indeed, by \ref{lem:cauchy-radius} we can choose three points $x_1, x_2, x_3\in \Sph(u,a)$ with $\lambda(x_i)$ from one $\Theta$-separated sequence
such that each distance $|w-x_i|$ is very close to~$b$.
Then choose $y,z\in \Sph(u,a)$ such that
\begin{align*}
|w-y|&<|w-x_i|,
&
|w-z|&<|w-x_i|
\end{align*}
for each $i$, and $|y-z|$ is almost $h(b)$.
Since $\Theta$ is a metric, we may choose $x=x_i$ that meets the conditions in \ref{eq:lambda-eps}.

Choose an almost shortest path $\gamma$ from $w$ to $x$, so its length is sufficiently close to
$|w-x|=b+\lambda(x)$.
Choose $b'>0$ (independently of $\eps$) such that $0\z<b-b'\z< \tfrac\delta{100}$; let $c'=a+b'$ and $w'\in\gamma\cap \Sph(u,c')$.
In particular, we have
\[
|w'-w|<\tfrac\delta{10}
\quad\text{and}\quad
|x-w'|<b'+2\cdot\lambda(x).
\]

\begin{wrapfigure}{o}{40mm}
\centering
\vskip-0mm
\includegraphics{mppics/pic-80}
\vskip-0mm
\end{wrapfigure}

Let $\rho$ be a reflection swapping $x$ with $y$;
note that $\rho(u)=u$.
Since
\begin{align*}
\tfrac\delta2&<\Theta(\lambda(x),\lambda(y))
\leqslant
\\
&\leqslant
|w-\rho(w)|=
\\
&=
2\cdot\inf\set{|w-o|}{o\in \Bis(x,y)},
\end{align*}
the point $w$ lies at distance at least $\tfrac\delta{10}$ from $\Bis(x,y)$.
Since $|w'-w|<\tfrac\delta{10}$, both points $w$ and $w'$ lie on the same side of $\Bis(x,y)$, and
\[|w-x|>|w-y|
\quad\Rightarrow\quad
|w'-x|>|w'-y|.\]
Therefore,
\[y\in \Sph(u,a)\cap \Ball(w',b'+2\cdot\lambda(x)),\]
and similarly,
\[z\in \Sph(u,a)\cap \Ball(w',b'+2\cdot\lambda(x)).\]
Thus \ref{sigma=const} follows.
\qeds

\subsection*{Expanding collinear sets}

\begin{thm}{Lemma}\label{cor:plus-one}
Let $f$ be a distance-preserving map from a finite set of real numbers $X$ into an all-set-homogeneous length space $M$.
Given three consecutive numbers $t_{i-1},t_i,t_{i+1}\in X$, set $t_i'=t_{i-1}-t_i +t_{i+1}$.
Then the map $f$ can be extended to a distance-preserving map on $X\cup\{t_i'\}$.
\end{thm}

\parit{Proof.}
Set $x_i=f(t_i)$ and let $f(t_i')\df x_i'\df\rho(x_i)$, where $\rho$ is a reflection that swaps $x_{i-1}$ and $x_{i+1}$.
We can assume $|x_{i-1}-x_i|\geqslant |x_i-x_{i+1}|$.
Then
\[|x_{i-1}-x'_i|+|x'_i-x_i|+|x_i-x_{i+1}|=|x_{i-1}-x_{i+1}|;\]
the triangle inequality gives the lower bound, and an almost shortest path from
$x_{i-1}$ to $x_i$ gives the upper bound.
Hence the statement follows.
\qeds

\begin{thm}{Proposition}\label{prop:ash>geodesic}
A complete all-set-homogeneous length space is geodesic.
\end{thm}

\parit{Proof.}
Let $M$ be a complete all-set-homogeneous length space.
Choose two points $u,v\in M$ and let $a=|u-v|$.
Assume $a<\diam M$.

Let $X$ be a finite set of real numbers $0=t_1<\ldots<t_n=a$,
and let $f\colon X\to M$ be a distance-preserving map such that $u=f(0)$ and $v=f(a)$.

By \ref{lem:plus-eps}, we can choose a point $w\in M$ such that $|u-v|+|v-w|=|u-w|$
and $b=|v-w|>0$ is arbitrarily small.
Let $c=a+b=|u-w|$ and let $f(c)=w$;
note that $f$ is still distance-preserving.

Applying \ref{cor:plus-one}, we can propagate a short gap from right to left,
extending $f$ to a larger subset of $[0,c]$ such that consecutive points are at distance at most~$b$.
Then we can restrict $f$ to the intersection of the enlarged domain with $[0,a]$ and repeat the procedure for a smaller $b>0$.

Repeating this procedure, we extend $f$ to a distance-preserving map on a dense set of $[0,a]$.
By completeness of $M$, the map $f$ extends to a geodesic from $u$ to $v$.
Thus a geodesic $[uv]$ exists if $|u-v|<\diam M$.

If $a=\diam M$, we can choose a sequence of geodesics $[uv_1],[uv_2],\ldots$ such that
$|u-v_n|\to a$ as $n\to\infty$.
By 2-point-homogeneity, we can assume that they are nested: $[uv_1]\subset[uv_2]\subset\ldots$
Taking their union and passing to its closure gives a geodesic of length $a$.
Now, by all-set-homogeneity, $M$ is geodesic.
\qeds

\section{Remarks}

In the thin-mirror case, one could prove geodesic existence and local uniqueness directly, using a family of pairwise commuting reflections of maximum cardinality and a version of the monster lemma adapted to the whole family.
The strengthened lemma would provide suitable short retractions onto half-intersections of mirrors.
The enlargement argument below then forces the relevant distance coordinate to be injective, reducing the construction of geodesics to an ordered length space.
We chose the present organization to make the geodesic case more accessible and independently readable.

\begin{thm}{Enlargement lemma}
Let $M$ be an all-set-homogeneous length space with thin mirrors.
Let $\rho_1,\ldots,\rho_n$ be distinct pairwise commuting reflections,
and $H_i=\Bis_{\rho_i}$.
Suppose distinct points $x,y\in H_2\cap\cdots\cap H_n$
lie on the same side and at the same distance from $H_1$.
Then $M$ admits $n+1$ distinct pairwise commuting reflections.
\end{thm}

\parit{Proof.}
Let $\sigma$ be the reflection swapping $x$ with $y$.
Note that $\sigma$ commutes with $\rho_2,\ldots,\rho_n$.

\begin{wrapfigure}{o}{30mm}
\centering
\vskip-4mm
\includegraphics{mppics/pic-100}
\vskip-0mm
\end{wrapfigure}

Set $u=\rho_1(x)$, $v=\rho_1(y)$, and $w=\sigma(u)$.
Note that $\rho_1$ commutes with $\sigma$ if and only if $v=w$;
in this case we can add $\sigma$ to the family, obtaining $n+1$ commuting reflections.

Thus we can assume that $v\ne w$; let $\tau$ be the reflection swapping $v$ with $w$.
The equality of the distances to $H_1$ gives $x,y\in\Bis(v,w)$.
Hence the reflection $\tau$ fixes $x,y$ and commutes with $\sigma$.
Each $\rho_i$, $i\geqslant2$, fixes $v$ and $w$; so it also commutes with $\tau$.
Thus, discarding $\rho_1$, we obtain $n+1$ commuting reflections
$\sigma$, $\tau$, $\rho_2,\ldots,\rho_n$.
\qeds

Classification of connected all-set-homogeneous spaces is the next natural step.
Such spaces appear if one modifies the metric on an all-set-homogeneous length space by a function.
Namely, let $M$ be a length space, set $L=\diam M$.
Suppose $\varphi\colon[0,L]\to[0,D]$ is a continuous strictly increasing subadditive function such that $\varphi(0)=0$; \emph{subadditive} means $\varphi(a+b)\leqslant \varphi(a)+\varphi(b)$ if both sides are defined.
Note that the function $(x,y)\mapsto \varphi(|x-y|)$ defines another metric on $M$;
let us denote by $M^\varphi$ the obtained metric space.
If $M$ is an all-set-homogeneous length space, then $M^\varphi$ is a connected all-set-homogeneous space.

\begin{thm}{Proposition}\label{prop:path-connected}
Any complete path-connected all-set-homogeneous space is isometric to $M^\varphi$, where $M$ is a universal metric tree of finite valence or a classical space and $\varphi\colon[0,\diam M]\to [0,D]$ is a continuous strictly increasing subadditive function such that $\varphi(0)=0$.
\end{thm}

\begin{thm}{Question}
Does the same hold for every connected all-set-homogeneous space?
\end{thm}

The proof of the proposition follows the idea of Jacques Tits \cite[IV.E § 3]{tits}.

\parit{Proof of \ref{prop:path-connected}.}
Let us show that any complete path-connected all-set-homogeneous space $X$ is isometric to $M^\varphi$, where $M$ is a length space and $\varphi$ is a strictly increasing function such that $\varphi(0)=0$.
Once this is done, then the main theorem finishes the proof.

Any ball (closed or open) in $X$ is path-connected.
Indeed, connect $p$ to $x$ by a path, say $\gamma$;
choose the first point $\gamma(t)$ at the distance $|p-x|$ from $p$;
choose an isometry $\iota$ that fixes $p$ and moves $\gamma(t)$ to $x$;
the path $\iota\circ\gamma|_{[0,t]}$ connects $p$ to $x$ and lies in the ball containing $x$ with the center at $p$.

Let $D=\diam X$.
Note that there is an operation $\oplus$ on $[0,D]$ such that
\[\bar\Ball(p,r\oplus s)=\Cl_X\left(\bigcup_{q\in \bar\Ball(p,r)}\bar\Ball(q,s)\right);\]
in particular, the set on the right hand side is a ball.
Observe that $\oplus$ is symmetric, associative, continuous, and $0<r<D$ $\Rightarrow$ $r<r\oplus r$.
By the linearization theorem \cite[(3.3)]{ling1965},
\[\varphi^{-1}(r\oplus s)=\min\{\,\varphi^{-1}(r)+\varphi^{-1}(s),L\,\}\]
for some increasing continuous subadditive function $\varphi\colon[0,L]\to [0,D]$.%
\footnote{We might have $L=\infty$ or $D=\infty$; for example, consider the function $\varphi\colon[0,\infty]\to[0,1]$ defined by $\varphi(t)=\tfrac t{t+1}$ with $\varphi(\infty)=1$.}
Note that
\[|x-z|\leqslant |x-y|\oplus |y-z|\leqslant|x-y|+ |y-z|\]
for any $x,y,z\in X$; therefore, subadditivity of $\varphi$ follows.

$\varphi^{-1}(|\ -\ |)$ defines a new metric on $X$; indeed any pair of points $x$ and $y$ can be connected by a chain $x=x_0,\dots,x_n=y$ such that $|x_{i}-x_{i-1}|<\eps$, where $\eps>0$ is small; therefore, $\varphi^{-1}(|x -y|)\leqslant n\cdot \varphi^{-1}(\eps)<\infty$; the remaining metric axioms are straightforward.
Denote by $M$ the resulting metric space; note that $L=\diam M$.
Since $X$ is complete, so is $M$.
Note that $M$ has almost midpoints; by Menger's lemma it is a length space.
\qeds

By Tits' theorem quoted in the introduction \cite{tits}, a locally compact geodesic 3-point-homogeneous space is classical; in particular it is  \emph{finite-set-homogeneous}; that is, it is $m$-point-homogeneous for every positive integer $m$.

The class of complete finite-set-homogeneous geodesic spaces is quite rich.
It includes infinite-dimensional analogs of classical spaces, as well as universal metric trees of infinite valence.
Another class of examples comes from metric Fraïssé limits;
in particular, Urysohn universal spaces (unbounded and bounded).
All these additional spaces are isometric to their own proper subsets (they also include infinite equilateral sets); therefore they are not all-set-homogeneous.

\begin{thm}{Question}
Is there a positive integer $n$ such that any $n$-point-homogeneous  geodesic space is finite-set-homogeneous?
\end{thm}

The champion so far is the space of Lebesgue-measurable sets up to measure zero on $[0,1]$ with metric defined by $\mu(A\triangle B)$, where $\triangle$ denotes the symmetric difference.
This space is 3-point-homogeneous, but not 4-point-homogeneous.
For general metric spaces, a negative answer was found by Manfred Droste and Dugald Macpherson \cite{droste-macpherson}.

\parbf{Acknowledgements.}
We thank Marcos Cossarini, James E. Hanson, Sergei Ivanov, Linus Kramer, and Alexander Lytchak for helpful discussions, both in person and on MathOverflow.
The proofs were developed in collaboration with AI (GPT-6 Astra).

{\sloppy
\def\emph{\textit}
\printbibliography[heading=bibintoc]\fussy}

@article {tits,
    AUTHOR = {Tits, J.},
     TITLE = {Sur certaines classes d'espaces homog\`enes de groupes de {L}ie},
   JOURNAL = {Acad. Roy. Belg. Cl. Sci. M\'{e}m. Coll. in 8$^\circ$},
  FJOURNAL = {Acad\'{e}mie Royale de Belgique. Classe des Sciences. M\'{e}moires.
              Collection in-8$^\circ$. Koninklijke Belgische Academie.
              Klasse der Wetenschappen. Verhandelingen. Verzameling
              in-8$^\circ$},
    VOLUME = {29},
      YEAR = {1955},
    NUMBER = {3},
     PAGES = {268},
      ISSN = {0365-0936},
   MRCLASS = {17.0X},
  MRNUMBER = {76286},
MRREVIEWER = {L. Auslander},
}

@misc{azuelos2025,
      title={Uncountably many homogeneous real trees with the same valence},
      author={P. Azuelos},
      year={2025},
      eprint={2511.03722},
      archivePrefix={arXiv},
      primaryClass={math.MG},
      url={https://arxiv.org/abs/2511.03722},
}

@article {bargetz-bartos-kubis-luggin,
AUTHOR = {Bargetz, C. and Barto\v s, A. and Kubi\'s, W. and Luggin, F.},
     TITLE = {Homogeneous isosceles-free spaces},
   JOURNAL = {Rev. R. Acad. Cienc. Exactas F\'is. Nat. Ser. A Mat. RACSAM},
  FJOURNAL = {Revista de la Real Academia de Ciencias Exactas, F\'isicas y
              Naturales. Serie A. Matematicas. RACSAM},
    VOLUME = {118},
      YEAR = {2024},
    NUMBER = {3},
     PAGES = {Paper No. 118, 32},
      ISSN = {1578-7303,1579-1505},
   MRCLASS = {54E35 (03C50 05C15 20B25)},
  MRNUMBER = {4748952},
MRREVIEWER = {Yasunao\ Hattori},
       DOI = {10.1007/s13398-024-01587-y},
       URL = {https://doi.org/10.1007/s13398-024-01587-y},
}

@article {birkhoff,
    AUTHOR = {Birkhoff, G.},
     TITLE = {Metric foundations of geometry. {I}},
   JOURNAL = {Trans. Amer. Math. Soc.},
  FJOURNAL = {Transactions of the American Mathematical Society},
    VOLUME = {55},
      YEAR = {1944},
     PAGES = {465--492},
      ISSN = {0002-9947},
   MRCLASS = {48.0X},
  MRNUMBER = {10393},
MRREVIEWER = {L. M. Blumenthal},
       DOI = {10.2307/1990304},
       URL = {https://doi.org/10.2307/1990304},
}

@book {busemann,
    AUTHOR = {Busemann, H.},
     TITLE = {Metric methods in {F}insler spaces and in the
              foundations of geometry},
    SERIES = {Annals of Mathematics Studies, No. 8},
 %PUBLISHER = {Princeton University Press, Princeton, N. J.},
      YEAR = {1942},
  %   PAGES = {viii+243},
   MRCLASS = {48.0X},
  MRNUMBER = {0007251},
MRREVIEWER = {S. M. Ulam},
}

@book {chiswell,
    AUTHOR = {Chiswell, I.},
     TITLE = {Introduction to {$\Lambda$}-trees},
% PUBLISHER = {World Scientific Publishing Co., Inc., River Edge, NJ},
      YEAR = {2001},
     %PAGES = {xii+315},
      ISBN = {981-02-4386-3},
   MRCLASS = {20E08 (03C60 20F65)},
  MRNUMBER = {1851337},
MRREVIEWER = {Vincent\ Guirardel},
       DOI = {10.1142/4495},
       URL = {https://doi.org/10.1142/4495},
}

@article {droste-macpherson,
    AUTHOR = {Droste, M. and Macpherson, H. D.},
     TITLE = {On {$k$}-homogeneous posets and graphs},
   JOURNAL = {J. Combin. Theory Ser. A},
  FJOURNAL = {Journal of Combinatorial Theory. Series A},
    VOLUME = {56},
      YEAR = {1991},
    NUMBER = {1},
     PAGES = {1--15},
      ISSN = {0097-3165,1096-0899},
   MRCLASS = {03C65 (03C35 03C50 05C25 05E25 06A06)},
  MRNUMBER = {1082839},
MRREVIEWER = {James\ H.\ Schmerl},
       DOI = {10.1016/0097-3165(91)90018-C},
       URL = {https://doi.org/10.1016/0097-3165(91)90018-C},
}

@article {dyubina-polterovich,
    AUTHOR = {Dyubina, A. and Polterovich, I.},
     TITLE = {Explicit constructions of universal {$\mathbb{R}$}-trees and
              asymptotic geometry of hyperbolic spaces},
   JOURNAL = {Bull. London Math. Soc.},
  FJOURNAL = {The Bulletin of the London Mathematical Society},
    VOLUME = {33},
      YEAR = {2001},
    NUMBER = {6},
     PAGES = {727--734},
      ISSN = {0024-6093},
   MRCLASS = {57M07 (20F67 53C23 54F50)},
  MRNUMBER = {1853785},
       DOI = {10.1112/S002460930100844X},
       URL = {https://doi.org/10.1112/S002460930100844X},
}

@article {lebedeva-petrunin2024,
    AUTHOR = {Lebedeva, N. and Petrunin, A.},
     TITLE = {All-set-homogeneous spaces},
   JOURNAL = {St. Petersburg Math. J.},
  FJOURNAL = {St. Petersburg Mathematical Journal},
    VOLUME = {35},
      YEAR = {2024},
    NUMBER = {3},
     PAGES = {473--476},
      ISSN = {1061-0022,1547-7371},
   MRCLASS = {54E35 (53C30)},
  MRNUMBER = {4907863},
MRREVIEWER = {Mikhail\ Ostrovskii},
       DOI = {10.1090/spmj/1814},
       URL = {https://doi.org/10.1090/spmj/1814},
}

@article {ling1965,
    AUTHOR = {Ling, C.-H.},
     TITLE = {Representation of associative functions},
   JOURNAL = {Publ. Math. Debrecen},
  FJOURNAL = {Publicationes Mathematicae Debrecen},
    VOLUME = {12},
      YEAR = {1965},
     PAGES = {189--212},
      ISSN = {0033-3883,2064-2849},
   MRCLASS = {39.30},
  MRNUMBER = {190575},
MRREVIEWER = {M.\ Hossz\'u},
       DOI = {10.5486/pmd.1965.12.1-4.19},
       URL = {https://doi.org/10.5486/pmd.1965.12.1-4.19},
}

@misc{niemiec2023,
      title={Extensive approach to absolute homogeneity},
      author={P. Niemiec},
      year={2023},
      eprint={2308.09986},
      archivePrefix={arXiv},
      primaryClass={math.GN},
      url={https://arxiv.org/abs/2308.09986},
}
\end{document}